\documentclass[11pt]{article}

\usepackage{amssymb,latexsym,amsmath}

\usepackage{graphicx}

\begin{document}

\newcommand{\bfi}{\bfseries\itshape}

\makeatletter

\@addtoreset{figure}{section}

\def\thefigure{\thesection.\@arabic\c@figure}

\def\fps@figure{h, t}

\@addtoreset{table}{bsection}

\def\thetable{\thesection.\@arabic\c@table}

\def\fps@table{h, t}

\@addtoreset{equation}{section}

\def\theequation{\thesubsection.\arabic{equation}}

\makeatother

\newtheorem{thm}{Theorem}[section]

\newtheorem{prop}[thm]{Proposition}

\newtheorem{lema}[thm]{Lemma}

\newtheorem{cor}[thm]{Corollary}

\newtheorem{defi}[thm]{Definition}

\newtheorem{rk}[thm]{Remark}

\newtheorem{exempl}{Example}[section]

\newenvironment{exemplu}{\begin{exempl}  \em}{\hfill $\surd$

\end{exempl}}

\newcommand{\comment}[1]{\par\noindent{\raggedright\texttt{#1}

\par\marginpar{\textsc{Comment}}}}

\newcommand{\todo}[1]{\vspace{5 mm}\par \noindent \marginpar{\textsc{ToDo}}\framebox{\begin{minipage}[c]{0.95 \textwidth}

\tt #1 \end{minipage}}\vspace{5 mm}\par}

\newcommand{\ea}{\mbox{{\bf a}}}

\newcommand{\eu}{\mbox{{\bf u}}}

\newcommand{\ueu}{\underline{\eu}}

\newcommand{\ueo}{\overline{u}}

\newcommand{\oeu}{\overline{\eu}}

\newcommand{\ew}{\mbox{{\bf w}}}

\newcommand{\ef}{\mbox{{\bf f}}}

\newcommand{\eF}{\mbox{{\bf F}}}

\newcommand{\eC}{\mbox{{\bf C}}}

\newcommand{\en}{\mbox{{\bf n}}}

\newcommand{\eT}{\mbox{{\bf T}}}

\newcommand{\eL}{\mbox{{\bf L}}}

\newcommand{\eR}{\mbox{{\bf R}}}

\newcommand{\eV}{\mbox{{\bf V}}}

\newcommand{\eU}{\mbox{{\bf U}}}

\newcommand{\ev}{\mbox{{\bf v}}}

\newcommand{\eve}{\mbox{{\bf e}}}

\newcommand{\uev}{\underline{\ev}}

\newcommand{\eY}{\mbox{{\bf Y}}}

\newcommand{\eK}{\mbox{{\bf K}}}

\newcommand{\eP}{\mbox{{\bf P}}}

\newcommand{\eS}{\mbox{{\bf S}}}

\newcommand{\eJ}{\mbox{{\bf J}}}

\newcommand{\eB}{\mbox{{\bf B}}}

\newcommand{\eH}{\mbox{{\bf H}}}

\newcommand{\leb}{\mathcal{ L}^{n}}

\newcommand{\eI}{\mathcal{ I}}

\newcommand{\eE}{\mathcal{ E}}

\newcommand{\hen}{\mathcal{H}^{n-1}}

\newcommand{\eBV}{\mbox{{\bf BV}}}

\newcommand{\eA}{\mbox{{\bf A}}}

\newcommand{\eSBV}{\mbox{{\bf SBV}}}

\newcommand{\eBD}{\mbox{{\bf BD}}}

\newcommand{\eSBD}{\mbox{{\bf SBD}}}

\newcommand{\ecs}{\mbox{{\bf X}}}

\newcommand{\eg}{\mbox{{\bf g}}}

\newcommand{\paromega}{\partial \Omega}

\newcommand{\gau}{\Gamma_{u}}

\newcommand{\gaf}{\Gamma_{f}}

\newcommand{\sig}{{\bf \sigma}}

\newcommand{\gac}{\Gamma_{\mbox{{\bf c}}}}

\newcommand{\deu}{\dot{\eu}}

\newcommand{\dueu}{\underline{\deu}}

\newcommand{\dev}{\dot{\ev}}

\newcommand{\duev}{\underline{\dev}}

\newcommand{\weak}{\stackrel{w}{\approx}}

\newcommand{\mild}{\stackrel{m}{\approx}}

\newcommand{\strong}{\stackrel{s}{\approx}}

\newcommand{\weakdown}{\rightharpoondown}

\newcommand{\opg}{\stackrel{\mathfrak{g}}{\cdot}}

\newcommand{\opunu}{\stackrel{1}{\cdot}}
\newcommand{\opdoi}{\stackrel{2}{\cdot}}

\newcommand{\opn}{\stackrel{\mathfrak{n}}{\cdot}}

\newcommand{\tr}{\ \mbox{tr}}

\newcommand{\Ad}{\ \mbox{Ad}}

\newcommand{\ad}{\ \mbox{ad}}

\renewcommand{\contentsname}{ }

\title{Construction  of bipotentials and a minimax theorem of Fan}

\author{Marius Buliga\footnote{"Simion Stoilow" Institute of Mathematics of the 
Romanian Academy, PO BOX 1-764,014700 Bucharest, Romania, e-mail: 
Marius.Buliga@imar.ro } ,
 G\'ery de Saxc\'e\footnote{Laboratoire de 
  M\'ecanique de Lille, UMR CNRS 8107, Universit\'e des Sciences et 
Technologies de Lille,
 Cit\'e Scientifique, F-59655 Villeneuve d'Ascq cedex, France, 
e-mail: gery.desaxce@univ-lille1.fr} , Claude  Vall\'ee\footnote{Laboratoire de 
M\'ecanique des Solides, UMR 6610, UFR SFA-SP2MI, bd M. et P. Curie, 
t\'el\'eport 2, BP 30179, 86962 Futuroscope-Chasseneuil cedex, 
France, e-mail: vallee@lms.univ-poitiers.fr} }

\date{This version: 19.02.2007}

\maketitle

\begin{abstract}
The bipotential theory is  based on an extension of Fenchel's inequality, 
with  several powerful applications related to non associated constitutive 
laws in Mechanics: frictional contact \cite{saxfeng}, non-associated Drucker-Prager 
model \cite{bersaxce},  or Lemaitre plastic ductile damage law \cite{bodo}, 
to cite a few. 
 
This is a second paper on the mathematics of the bipotentials, 
following \cite{bipo1}. We prove here another  reconstruction theorem for a 
bipotential  from a convex lagrangian cover, this time using a  
convexity notion related to   a minimax theorem of Fan. 
\end{abstract}


{\bf MSC-class:} 49J53; 49J52; 26B25

\section{Introduction}

In Mechanics, the theory of standard materials is a well-known application 
of Convex Analysis. However, the so-called  non-associated constitutive laws  
cannot be cast in the mould of the standard materials. 

From the mathematical viewpoint, a  non associated constitutive law 
is  a multivalued operator $\displaystyle T:X\rightarrow 2^{Y}$ 
{\it which is not supposed to be monotone}. Here  $X$, $Y$ are dual locally 
convex spaces, with duality product $\langle \cdot , \cdot \rangle: X \times Y \rightarrow \mathbb{R}$.

A possible way to study  non-associated constitutive laws by using 
 Convex Analysis,   proposed first in \cite{saxfeng}, 
consists in constructing a "bipotential" function $b$ of two variables, 
which physically represents the dissipation.  

A bipotential  function $b$  is 
 bi-convex, satisfies an inequality generalizing Fenchel's one, 
$ \forall x \in X, y \in Y, \: \ b(x,y) \geq \langle x, y \rangle $, and a 
relation involving partial subdifferentials of $b$ with respect to  variables  
$x$, $y$.     
In the case of associated constitutive laws  the 
bipotential has the expression $b(x,y) = \phi(x)+\phi^{*}(y)$. 

The graph of a bipotential $b$ is simply the set $M(b) \subset X \times Y$ of those 
pairs $(x,y)$ such that $ b(x,y) =  \langle x, y \rangle $. A 
multivalued operator  $\displaystyle T:X\rightarrow 2^{Y}$ is expressed with 
the help of the bipotential $b$ if the graph of $T$ (in the usual sense) equals 
$M(b)$. 

The  non associated constitutive laws which can be expressed with the 
help of bipotentials are called in Mechanics implicit, or weak, normality 
rules. They  have the  form of an implicit relation between dual 
variables,  $ y \in \partial b(\cdot , y)(x) $. 
 
Among the  applications of bipotentials to Solid Mechanics we cite: 
Coulomb's friction law \cite{sax CRAS 92} , non-associated Dr\"ucker-Prager 
\cite{sax boush KIELCE 93}  and Cam-Clay models \cite{sax BOSTON 95}  in 
Soil Mechanics, cyclic Plasticity (\cite{sax CRAS 92},\cite{bodo sax EJM 01}) 
and Viscoplasticity \cite{hjiaj bodo CRAS 00} of metals with non linear 
kinematical hardening rule, Lemaitre's damage law \cite{bodo}, the coaxial 
laws (\cite{dangsax},\cite{vall leri CONST 05}). A  review of these laws 
expressed in terms of bipotentials  can be found in \cite{dangsax} and 
\cite{vall leri CONST 05}.

In order to better understand the bipotential approach, 
in the paper \cite{bipo1} we solved two key problems: (a) when the graph of a 
given multivalued operator can be expressed as the set of critical points of a 
bipotentials, and (b) a method of construction of a bipotential associated 
(in the sense of point (a)) to a multivalued, typically non monotone, operator. 

Our main tool was the notion of convex lagrangian cover of the graph of 
the multivalued operator, and a related notion of implicit convexity of this 
cover. 

In this paper we prove another  reconstruction theorem for a 
bipotential  from a convex lagrangian cover, this time using a  
convexity notion related to   a minimax theorem of Fan.

\section{Notations and definitions}

$X$ and $Y$ are topological, locally convex, real vector spaces of dual 
variables $x \in X$ and $y \in Y$, with the duality product 
$\langle \cdot , \cdot \rangle : X \times Y \rightarrow \mathbb{R}$. 
The topologies of the spaces $X, Y$ are  compatible with the duality 
product, that is: any  continuous linear functional on $X$ (resp. $Y$) 
has the form $x \mapsto \langle x,y\rangle$, for some $y \in Y$ (resp. 
$y \mapsto \langle x,y\rangle$, for some  $x \in X$).

We use the notation: $\displaystyle \bar{\mathbb{R}} = \mathbb{R}\cup \left\{ +\infty \right\}$. 

Given a function 
$\displaystyle \phi: X \rightarrow  \bar{\mathbb{R}}$, the domain 
$dom \, \phi$ is the set of points with value other than $+\infty$. 
The polar of $\phi$, or Fenchel conjugate, 
$\phi^{*}: Y \rightarrow \bar{\mathbb{R}}$ 
is defined by: 
$\displaystyle \phi^{*}(y) = \sup \left\{ \langle y,x\rangle - 
\phi(x) \mid x \in X \right\}$.  

We denote by $\Gamma(X)$ the class of convex and lower semicontinuous 
functions $\displaystyle \phi: X \rightarrow \bar{\mathbb{R}}$. 
The  class of convex and lower semicontinuous functions 
$\displaystyle \phi: X \rightarrow \mathbb{R}$ is denoted by 
$\displaystyle \Gamma_{0}(X)$.

The subdifferential   of a function 
$\displaystyle \phi: X \rightarrow \bar{\mathbb{R}}$ in a point 
$x \in dom \, \phi$ is the (possibly empty) set: 
$$\partial \phi(x) = \left\{ u \in Y \mid \forall z \in X  \  
\langle z-x, u \rangle \leq \phi(z) - \phi(x) \right\} \  .$$ 
In a similar way is defined the subdifferential  of a function 
$\psi: Y \rightarrow \bar{\mathbb{R}}$ in a point $y \in dom \, \psi$, 
as the set: 
$$\partial \psi(y) = \left\{ v \in X \mid \forall w \in Y  \  \langle v, w-y \rangle \leq \psi(w) - \psi(y) \right\} \ .$$ 

With these notations we have the Fenchel inequality: let 
$\displaystyle \phi: X \rightarrow \bar{\mathbb{R}}$ be a convex lower 
semicontinuous function. Then: 
\begin{enumerate}
\item[(i)] for any $x \in X , y\in Y$ we have $\displaystyle \phi(x) + \phi^{*}(y)  \geq \langle x, y \rangle$; 
\item[(ii)]  for any $(x,y) \in X \times Y$ we have the equivalences: 
$$ y \in \partial \phi(x) \ \Longleftrightarrow \ x \in \partial \phi^{*}(y)  \ \Longleftrightarrow \  \phi(x) + \phi^{*}(y)  = 
\langle x , y \rangle \ . $$
\end{enumerate}

\begin{defi} To a graph  $M \subset X \times Y$ we associate the multivalued operators: 

$$\displaystyle  X \ni x  \mapsto m(x) \ = \ \left\{ y \in Y \mid (x,y) \in 
M \right\} \ ,$$
$$\displaystyle  Y \ni y  \mapsto m^{*}(y) \ = \ \left\{ x \in X \mid (x,y) \in M \right\} \ .$$
The domain of the graph $M$  is by definition  $\displaystyle dom(M) = \left\{ x \in 
X \mid m(x) \not = \emptyset\right\}$. 
The image  of the graph  $M$  is the set $\displaystyle im(M) = \left\{ 
y \in Y \mid m^{*}(y) \not = \emptyset\right\}$. 
\label{def1}
\end{defi}

\section{Bipotentials}

The notions and results  in this section were introduced  or 
proved in \cite{bipo1}. 

\begin{defi} A bipotential is a function $b: X \times Y \rightarrow 
\bar{\mathbb{R}}$ with the properties: 
\begin{enumerate}
\item[(a)] $b$ is convex and lower semicontinuos in each argument; 
\item[(b)] for any $x \in X , y\in Y$ we have $\displaystyle b(x,y) \geq 
\langle x, y \rangle$; 
\item[(c)]  for any $(x,y) \in X \times Y$ we have the equivalences: 
\begin{equation}
y \in \partial b(\cdot , y)(x) \ \Longleftrightarrow \ x \in \partial 
b(x, \cdot)(y)  \ \Longleftrightarrow \ b(x,y) = 
\langle x , y \rangle \ .
\label{equiva}
\end{equation}
\end{enumerate}
The graph of $b$ is 
\begin{equation}
M(b) \ = \ \left\{ (x,y) \in X \times Y \ \mid \ b(x,y) = \langle x, y 
\rangle \right\} \  .
\label{mb}
\end{equation}
\label{def2}
\end{defi}

{\bf Examples.} {\bf (1.)} (Separable bipotential) If 
$\phi: X \rightarrow \mathbb{R}$ is a convex, lower semicontinuous  potential, 
consider the multivalued operator $\displaystyle \partial \phi$ (the 
subdifferential of $\phi$).  The graph of this operator  is the set 
\begin{equation}
M(\phi) \ = \ \left\{ (x,y) \in X \times Y \ \mid \ \phi(x)+\phi^{*}(y) = \langle x, y \rangle \right\} \  .
\label{mphi}
\end{equation}
$M(\phi)$ is {\it maximally  cyclically 
monotone} \cite{rocka} Theorem 24.8. Conversely, if 
$M$ is closed and maximally cyclically monotone then there is a convex, 
lower semicontinuous $\phi$ such that $M=M(\phi)$.  

To  the function  $\phi$  we  associate the 
{\it separable bipotential} $$\displaystyle b(x,y)  = \phi(x) + \phi^{*}(y) .$$ 
Indeed, the Fenchel inequality can be reformulated by saying that the function 
$b$,  previously defined, is a bipotential. More precisely, the point (b) 
(resp. (c))  in the definition of a bipotential corresponds to (i) 
(resp. (ii)) from Fenchel inequality.

The bipotential $b$ and the function  
$\phi$ have the same graph: $\displaystyle M(b) = M(\phi)$. 

{\bf (2.)} (Cauchy bipotential)   Let $X=Y$ be a Hilbert space and let the duality 
product be  equal 
to the scalar product. Then we define the {\it Cauchy bipotential} by the formula 
$$\displaystyle b(x,y) = \| x\| \  \| y\| .$$ 
Let us check the Definition (\ref{def2}) The point (a) is obviously satisfied. 
The point (b) is true by the Cauchy-Schwarz-Bunyakovsky inequality.  
We have equality in the Cauchy-Schwarz-Bunyakovsky inequality 
$b(x,y) = \langle x,y \rangle$ if and only if there is $\lambda > 0$ such 
that $y = \lambda x$ or one of $x$ and $y$ vanishes.  This is exactly the statement from the  point (c), for 
the function $b$ under study.

 The graph $M(b)$ is the set of pairs of collinear 
and with same orientation vectors. It can not be expressed by a separable 
bipotential because $M(b)$ is not a cyclically monotone graph. 

\begin{defi}  The non empty set $M \subset X \times Y$ is a  BB-graph  
(bi-convex, bi-closed) if for all $x \in \ dom(M)$ and for all $y \in \ im(M)$ 
the sets $\displaystyle m(x)$ and $\displaystyle m^{*}(y)$ are convex and closed.
\label{dh1}
\end{defi}

The following theorem  gives a necessary and   sufficient condition for the 
existence of a bipotential associated to a constitutive law $M$. 

\begin{thm}
 Given  a non empty set $M \subset X \times Y$, there is a bipotential 
$b$ such that $M=M(b)$ if and only if $M$ is a BB-graph. 
 \label{thm1}
 \end{thm}

Given the BB-graph graph $M$, the uniqueness of bipotential $b$ such that 
$M=M(b)$ {\it is not true}. For example, in the case of the Cauchy bipotential 
$b$, the proof of theorem \ref{thm1} (theorem ... \cite{bipo1}) 
provides a bipotential, denoted by $\displaystyle b_{\infty}$, such that 
$\displaystyle M(b) = M(b_{\infty})$ but $\displaystyle b \not = b_{\infty}$. 
This is 
in contrast with the case of a maximal cyclically monotone graph $M$, when 
by Rockafellar theorem (\cite{rocka} Theorem 24.8.) we have a method to 
reconstruct unambigously the associated separable bipotential. 

We noticed that in mechanical applications, we were able to 
reconstruct the physically relevant bipotentials $b$ starting from $M(b)$, by 
knowing a little more than the graph $M(b)$. This supplementary information is encoded in the following notion.

\begin{defi} Let $M \subset X \times Y$ be a non empty set.  A 
 convex lagrangian cover of  $M$ is a function   
$\displaystyle \lambda \in \Lambda \mapsto \phi_{\lambda}$ from  $\Lambda$ with 
values in the set  $\Gamma(X)$, with the 
properties:
\begin{enumerate}
\item[(a)] The set $\Lambda$ is a non empty compact topological space, 
\item[(b)] Let $f: \Lambda \times X \times Y \rightarrow \bar{\mathbb{R}}$ be the function defined by 
$$f(\lambda, x, y) \ = \ \phi_{\lambda}(x) + \phi^{*}_{\lambda}(y) .$$
Then for any $x \in X$ and for any $y \in Y$ the functions 
$f(\cdot, x, \cdot): \Lambda \times Y \rightarrow \bar{\mathbb{R}}$ and 
$f(\cdot, \cdot , y): \Lambda \times X \rightarrow \bar{\mathbb{R}}$ are  lower 
semi continuous  on the product spaces   $\Lambda \times Y$ and respectively 
$\Lambda \times X$ endowed with the standard topology, 
\item[(c)] We have 
$$M  \ = \  \bigcup_{\lambda \in \Lambda} M(\phi_{\lambda}) \quad  .$$
\end{enumerate}
\label{defcover}
\end{defi}

Not 
any BB-graph admits a convex lagrangian cover. There exist   BB-graphs  admitting only 
one  convex lagrangian cover (up to reparametrization), as well as BB-graphs 
which have infinitely many lagrangian covers.  The problem of describing 
the set of all convex lagrangian covers of a BB-graph seems to be difficult. 
We shall not discuss this problem here, but see the sections 5 and 8 
 in \cite{bipo1}.

The  
results in this paper  apply only to BB-graphs admitting at least one convex 
lagrangian cover.

 To a  a convex lagrangian cover we associate a function which will turn out to 
be a bipotential, under some supplementary hypothesis.

\begin{defi}
Let $\displaystyle \lambda \mapsto \phi_{\lambda}$ be a convex lagrangian 
cover of the BB-graph $M$. To the cover we associate the function
$b: X \times Y \rightarrow \mathbb{R} \cup \left\{ + \infty \right\}$ 
by the formula 
$$b(x,y) \ = \ \inf \left\{ \phi_{\lambda}(x)+ \phi_{\lambda}^{*}(y) 
\mbox{ : } \lambda \in \Lambda\right\} \ = \  \inf_{\lambda \in \Lambda}  
f(\lambda, x, y)  \quad  . $$
\label{defrecipe}
\end{defi}

In \cite{bipo1} we  imposed an implicit convexity inequality in order to get 
a function $b$ which is a bipotential. We need two definitions.

\begin{defi}
Let $\Lambda$ be an arbitrary non empty set and $V$ a real vector space. The 
function $f:\Lambda\times V \rightarrow \bar{\mathbb{R}}$ is 
implicitly  convex if for any two elements 
$\displaystyle (\lambda_{1}, z_{1}) , 
(\lambda_{2},  z_{2}) \in \Lambda \times V$ and for any two numbers 
$\alpha, \beta \in [0,1]$ with $\alpha + \beta = 1$ there exists 
$\lambda  \in \Lambda$ such that 
$$f(\lambda, \alpha z_{1} + \beta z_{2}) \ \leq \ \alpha 
f(\lambda_{1}, z_{1}) + \beta f(\lambda_{2}, z_{2}) \quad .$$
\label{defimpl}
\end{defi}

\begin{defi}
Let $\displaystyle \lambda \mapsto \phi_{\lambda}$ be a convex lagrangian cover 
of the BB-graph $M$ and $f: \Lambda \times X \times Y \rightarrow \mathbb{R}$ the associated function introduced in Definition \ref{defcover}, that is the function defined by 
$$f(\lambda, z, y) \ = \ \phi_{\lambda}(z) + \phi^{*}_{\lambda}(y) \quad  .$$ 
The cover is bi-implicitly convex (or a   BIC-cover) if 
for any $y \in  Y$ and $x\in X$  the functions $f(\cdot, \cdot, y)$ and 
$f(\cdot, x, \cdot)$ are implicitly convex in the sense of Definition 
\ref{defimpl}. 
\label{defbit}
\end{defi}

In the case of $M=M(\phi)$, with $\phi$ convex and lower semi continuous  
(this corresponds to separable bipotentials), the set $\Lambda$ has only one 
element $\Lambda = \left\{ \lambda \right\}$ 
  and we have only one potential $\displaystyle \phi$. The associated 
bipotential from Definition \ref{defrecipe} is obviously 
$$b(x,y) \ = \ \phi(x) + \phi^{*}(y) \ .$$
This is a BIC-cover in a trivial way: the implicit convexity conditions are equivalent with 
the convexity of $\displaystyle \phi$, $\displaystyle \phi^{*}$ respectively.

With this convexity condition we obtained in \cite{bipo1} the following result. 

\begin{thm} Let $\displaystyle \lambda \mapsto \phi_{\lambda}$ be a BIC-cover of 
 the BB-graph $M$ and $b: X \times Y \rightarrow R$ defined by
\begin{equation}
b(x,y) \ = \ \inf \left\{ \phi_{\lambda}(x) + \phi^{*}_{\lambda}(y) \ \mid \  \lambda \in \Lambda \right\} \ . 
\end{equation}
Then $b$ is a bipotential and $M=M(b)$. 
\label{thm2}
\end{thm}

\section{Main result}

For simplicity, in this section we shall work only with lower semi continuous convex functions 
$\phi$ with the property that $\phi$ and its Fenchel dual 
$\displaystyle \phi^{*}$  
take values 
in $\mathbb{R}$. 

We reproduce here 
the following definition of convexity (in a generalized sense), given by 
K. Fan \cite{kyfan} p. 42.

\begin{defi} Let $X$, $Y$ be two  arbitrary non empty sets. The function   
$f:X\times Y  \rightarrow \mathbb{R}$ is convex on $X$ in the sense of Fan if for any two 
elements $\displaystyle x_{1}, x_{2} \in X$ and for any two numbers 
$\alpha, \beta \in [0,1]$ with $\alpha + \beta = 1$ there exists a $x \in X$ 
such that for all $y \in Y$: 
$$f(x,y) \ \leq \ \alpha f(x_{1},y) + \beta f(x_{2},y) .$$
\label{defkyfan}
\end{defi}

With the help of the previous definition we introduce a new convexity condition 
for a convex lagrangian cover.

\begin{defi}
Let $\displaystyle \lambda \mapsto \phi_{\lambda}$ be a convex lagrangian cover 
of the BB-graph $M$. Consider the 
functions: 
$$g: X \times \Lambda  \times X \rightarrow \mathbb{R} \quad , \quad 
h: Y \times \Lambda \times Y \quad , $$
given by $\displaystyle g(x,\lambda, z) = \phi_{\lambda}(x) - \phi_{\lambda}(z)$, 
respectively  $\displaystyle h(y, \lambda, u) = \phi_{\lambda}^{*}(y) - 
\phi_{\lambda}^{*}(u)$. 

The cover is Fan bi-implicitly convex  if for any $x \in X$, $y \in Y$, 
the functions $g(x, \cdot , \cdot)$, $h(y, \cdot , \cdot)$ are convex in the 
sense of Fan on $\Lambda \times X$, $\Lambda \times Y$ respectively. 
\label{deffic}
\end{defi}

Recall the following  minimax theorem of Fan \cite{kyfan}, Theorem 2. 
In the formulation of the theorem words "convex" and "concave" have the meaning given in definition \ref{defkyfan} 
(more precisely $f$ is concave if $-f$ is convex in the sense of the before 
mentioned definition). 

\begin{thm} (Fan) Let $X$ be a compact Hausdorff space and $Y$ an arbitrary 
set. Let $f$ be a real valued function on $X\times Y$ such that, for every $y 
\in Y$, $f(\cdot ,y)$ is lower semicontinuous on $X$. If $f$ is convex on $X$ 
and concave on $Y$, then the expressions 
$\displaystyle \min_{x \in X} \sup_{y \in Y} f(x,y)$ and $\displaystyle 
\sup_{y \in Y} \min_{x \in X} f(x,y)$ have meaning, and 
$$\min_{x \in X} \sup_{y \in Y} f(x,y) = \sup_{y \in Y} \min_{x \in X} f(x,y) 
\quad . $$
\label{tfan}
\end{thm}

The difficulty of theorem \ref{thm2} boils down to the fact the  class of 
convex functions is not closed with respect to the inf operator. Nevertheless, 
by using Fan theorem \ref{tfan} we get the following general result.

\begin{thm} Let $\Lambda$ be a compact Hausdorff space and 
$\displaystyle \lambda \mapsto \phi_{\lambda} \in \Gamma_{0}(X)$ be a 
convex lagrangian cover  of the BB-graph $M$ such that: 
\begin{enumerate}
\item[(a)]   for any $x \in X$ and for any $y \in Y$ the functions 
$\displaystyle \Lambda \ni \lambda  \mapsto \phi_{\lambda}(x) \in \mathbb{R} $ 
and $\displaystyle \Lambda \ni \lambda  \mapsto 
\phi_{\lambda}^{*}(y) \in \mathbb{R}$  are  continuous,  
\item[(b)] the cover is Fan bi-implicitely convex in the sense of definition 
\ref{deffic}.
\end{enumerate}
Then  the function  
$b: X \times Y \rightarrow \mathbb{R}$ defined by
$$b(x,y)  \ = \ \inf \left\{ \phi_{\lambda}(x) + \phi^{*}_{\lambda}(y) \ \mid \  \lambda \in \Lambda \right\}$$
is a bipotential and $M=M(b)$. 
\label{thm4}
\end{thm}

\paragraph{Proof.} 
For some of the details of the proof we refer to the proof of theorem
 \ref{thm2}    in  \cite{bipo1} (in that paper theorem 4.12). 
There are five steps in that proof. In order to prove our theorem we have only
 to modify the first two steps: we want to show that for any $x \in  dom(M)$ 
and any $y \in  im(M)$ the functions  $b(\cdot, y)$ and $b(x, \cdot)$ are 
convex and lower semi continuous.
 
For $(x,y) \in X \times Y$ let us define the function 
$\displaystyle \overline{xy}: \Lambda \times X \rightarrow \mathbb{R}$ by 
$$\overline{xy} (\lambda, z) = \langle z,y\rangle + \phi_{\lambda}(x) - \phi_{\lambda}(z)  \quad . $$
We check now that $\displaystyle \overline{xy}$ verifies the hypothesis of 
theorem \ref{tfan}. Indeed, the hypothesis (a) implies that  for any $z \in X$ 
the function  $\displaystyle \overline{xy}(\cdot , z)$  is 
continuous. Notice that 
$$ \overline{xy} (\lambda,z) = \langle z, y \rangle + g(x,\lambda,z) \quad . $$
It follows from hypothesis (b) that the function $\displaystyle \overline{xy}$ is convex on 
$\Lambda$ in the sense of Fan. 

In order to prove the concavity of  
$\displaystyle \overline{xy}$ on $X$, it 
suffices to show that for any $\displaystyle z_{1}, z_{2} \in X$, for any 
$\alpha, \beta \in [0,1]$ such that $\alpha + \beta = 1$, we have the inequality 
$$\overline{xy}(\lambda, \alpha z_{1} + \beta z_{2}) \leq 
\alpha \overline{xy}(\lambda, z_{1}) + \beta \overline{xy}(\lambda, z_{2}) $$ 
for any $\lambda \in \Lambda$.   This inequality  is equivalent with 
$$\langle \alpha z_{1} + \beta z_{2} , y \rangle - \phi_{\lambda}(\alpha z_{1} + 
\beta z_{2}) \leq \alpha \left( \langle z_{1} , y\rangle - \phi_{\lambda}(z_{1} \right) + 
\beta \left( \langle z_{2} , y\rangle - \phi_{\lambda}(z_{2} \right)$$
for any $\lambda \in \Lambda$. But this  is implied by  the 
convexity of $\displaystyle \phi_{\lambda}$ for any $\lambda \in \Lambda$. 

In conclusion the function $\displaystyle \overline{xy}$ satisfies the hypothesis 
of theorem \ref{tfan}.  We deduce that 
$$\min_{\lambda  \in \Lambda} \sup_{z \in X}  \overline{xy}(\lambda, z) = \sup_{z \in X} \min_{\lambda \in \Lambda}\overline{xy}(\lambda, z) \quad . $$
Let us compute  the two sides of this equality. 

For the left hand side (LHS) we have: 
$$LHS = min_{\lambda  \in \Lambda} \sup_{z \in X}  \left\{ \langle z,y\rangle + \phi_{\lambda}(x) - \phi_{\lambda}(z) \right\} = $$
$$= \min_{\lambda  \in \Lambda}  \left\{ \phi_{\lambda}(x) + \sup_{z \in X} \left\{ \langle z,y\rangle - \phi_{\lambda}(z) \right\}  \right\} = $$
$$ =  \min_{\lambda  \in \Lambda}  \left\{ \phi_{\lambda}(x) + \phi^{*}_{\lambda}(y) \right\} = b(x,y) \quad . $$
For the right hand side (RHS) we have: 
$$RHS = \sup_{z \in X} \min_{\lambda \in \Lambda}\left\{ \langle z,y\rangle + \phi_{\lambda}(x) - \phi_{\lambda}(z) \right\} = $$
$$= \sup_{z \in X}  \left\{  \langle z,y\rangle -  \max_{\lambda \in \Lambda}\left\{  \phi_{\lambda}(z) - \phi_{\lambda}(x) \right\} \right\} \quad . $$
Let $\displaystyle \overline{x}: X \rightarrow \mathbb{R}$ be the function 
$$ \overline{x}(z) = \max_{\lambda \in \Lambda}\left\{  \phi_{\lambda}(z) - \phi_{\lambda}(x) \right\} \quad . $$
Then the right hand side RHS is in fact: 
$$RHS = \overline{x}^{*}(y) \quad . $$
Therefore we proved the equality: 
$$b(x,y) =   \overline{x}^{*}(y) \quad . $$
This shows that the function $b$ is convex and lower semicontinuous in the second argument. 

In order to 
prove that $b$ is convex and lower semicontinuous in the first  argument,  
replace $\displaystyle \phi_{\lambda}$ by $\displaystyle 
\phi_{\lambda}^{*}$ in the previous reasoning. $\quad \blacksquare$

\vspace{.5cm}

\vspace{.5cm}

\vspace{\baselineskip}

\end{document}